\documentclass[]{imsart}

\usepackage[]{graphicx}\usepackage[]{color}

\makeatletter
\def\maxwidth{ %
  \ifdim\Gin@nat@width>\linewidth
    \linewidth
  \else
    \Gin@nat@width
  \fi
}
\makeatother

\definecolor{fgcolor}{rgb}{0.345, 0.345, 0.345}

\usepackage{framed}
\makeatletter
 {\par\unskip\endMakeFramed%
 \at@end@of@kframe}
\makeatother

\definecolor{shadecolor}{rgb}{.97, .97, .97}
\definecolor{messagecolor}{rgb}{0, 0, 0}
\definecolor{warningcolor}{rgb}{1, 0, 1}
\definecolor{errorcolor}{rgb}{1, 0, 0}

\usepackage{alltt} 
\usepackage[utf8]{inputenc}
\usepackage{enumitem}
\usepackage{todonotes}
\usepackage{latexsym}
\usepackage{amssymb}
\usepackage{epsfig}
\usepackage{amsmath}
\usepackage{xcolor}
\usepackage{float}
\usepackage{hyperref}
\usepackage{subfig}
\usepackage{multirow}
\usepackage{booktabs}
\usepackage[title]{appendix}

\IfFileExists{upquote.sty}{\usepackage{upquote}}{}

\usepackage{enumitem}
\usepackage{todonotes}
\RequirePackage[OT1]{fontenc}
\RequirePackage{amsmath,amssymb,amsthm}
\RequirePackage[round]{natbib}

\hypersetup{
  citecolor=black}
	
\RequirePackage[colorlinks,citecolor=blue,urlcolor=blue]{hyperref}
\usepackage{graphicx}
\usepackage{mathtools}

\startlocaldefs
\numberwithin{equation}{section}
\theoremstyle{plain}

\endlocaldefs

\theoremstyle{plain}
\long\def\comment#1{}

\newtheorem{theorem}{Theorem}

\theoremstyle{definition}
\newtheorem{definition}{Definition}
\numberwithin{definition}{section}

\numberwithin{remark}{section}

\newcommand{\R}{\mathbb{R}}

\renewcommand{\qed}{\hfill{\tiny \ensuremath{\blacksquare} }}%

\newcommand{\E}{\mathbb{E}}

\makeatletter
\newcommand{\setword}[2]{%
  \phantomsection
  #1\def\@currentlabel{\unexpanded{#1}}\label{#2}%
}
\makeatother

\begin{document}

\begin{frontmatter}
\title{A Note on High-Dimensional Confidence Regions}
\runtitle{Volume of High-Dimensional Confidence Regions}
\runauthor{Klaassen}
\thankstext{T1}{Version August 2020.}

\begin{aug}
\author{\fnms{\textit{Sven}} \snm{\textit{Klaassen}}\ead[label=e2]{}\thanksref{T2}},

\thankstext{T2}{Corresponding author: \textit{sven.klaassen@uni-hamburg.de}.}

\affiliation{University of Hamburg\thanksmark{m1}}

\address{University of Hamburg}

\end{aug}

\begin{abstract}
Recent advances in statistics introduced versions of the central limit theorem for high-dimensional vectors, allowing for the construction of confidence regions for high-dimensional parameters. In this note, $s$-sparsely convex high-dimensional confidence regions are compared with respect to their volume. Specific confidence regions which are based on $\ell_p$-balls are found to have exponentially smaller volume than the corresponding hypercube. The theoretical results are validated by a comprehensive simulation study.
\end{abstract}


\end{frontmatter}

\section{Introduction}
Constructing valid confidence regions is essential to assess the uncertainty which is associated with point estimates. Even for a single parameter there exist multiple valid confidence intervals. As a result, a large body of literature has been developed to construct confidence intervals which have desirable properties. In general, confidence intervals are constructed to minimize corresponding volume (see for example \citet{efron2006minimum} and \citet{jeyaratnam1985minimum}). 
In recent advances, \citet{chernozhukov2017central} developed a central limit theorem for a high-dimensional vector of random variables, allowing for confidence regions for high-dimensional parameter vectors. Nevertheless, these results only hold for specific sets $\mathcal{A}$ which are not too complex. A common application (see e.g. \citet{belloni2018uniformly}) of their work is to construct confidence regions in shape of a hypercube. 
We consider a more general setting of $s$-sparsely convex confidence regions which are then shown to have exponentially smaller volume than the corresponding cube.
\section{Notation}
Throughout the paper, we consider a random element $X$ from some common probability space $(\Omega,\mathcal{A},P)$.
We denote by $P\in \mathcal{P}_n$ a probability measure out of large class of probability measures, which may vary with the sample size (since the model is allowed to change with $n$) and by $\mathbb{P}_n$ the empirical probability measure. Additionally, let $\E$, respectively $\E_n$, be the expectation with respect to $P$, respectively $\mathbb{P}_n$.
For an $\alpha\in (0;\frac{1}{2})$ and a real valued random variable $Z$, define
\begin{align*}
q_\alpha(Z):= (1-\alpha)\text{ - Quantile of } Z.
\end{align*}
For a given set $A\in\R^d$, define
\begin{align*}
V(A):=\int\dots\int 1_{A} dx_1\dots dx_d.
\end{align*}
 Further, for a vector $v\in \mathbb{R}^d$ and $p\ge 1$ denote the $\ell_p$ norm 
\begin{align*}
\|v\|_p:=\left(\sum_{l=1}^p |v_l|^p\right)^{1/p},
\end{align*} 
$\|v\|_0$ equals the number of non-zero components and $\|v\|_\infty=\sup_{l=1,\dots,d }|v_l|$ denotes the $\sup$-norm. For any subset $J=\{J_1,\dots,J_k\}\subseteq \{1,\dots,d\}$, we define
$$v_{J}:=(v_{J_1},\dots,v_{J_{k}})^T\in\R^{k}$$ as the corresponding subvector of $v$.\\
Let $A$ be a $m\times d$ matrix. Denote the operator norms on $\mathbb{R}^{m\times d}$, which are induced by the $\ell_p$ vector norms, $\|A\|_{q,p} := \sup_{v\in \mathbb{R}^d:\|v\|_q= 1} \|Av\|_p$ and $\|A\|_{p} := \|A\|_{p,p}$.\\
Let $c$ and $C$ denote positive constants independent of $n$ with values that may change at each appearance. The notation $a_n\lesssim b_n$ means $a_n\le Cb_n$ for all $n$ and some $C$. Furthermore, $a_n=o(1)$ denotes that there exists a sequence $(b_n)_{\ge 1}$ of positive numbers such that $|a_n|\le b_n$ for all $n$ where $b_n$ is independent of $P\in\mathcal{P}_n$ for all $n$ and $b_n$ converges to zero. Finally, $a_n=O_P(b_n)$ means that for any $\epsilon>0$, there exists a $C$ such that $P(a_n>Cb_n)\le\epsilon$ for all $n$.
\section{Main Results}
At first, we recap the high-dimensional limit theorem from \citet{chernozhukov2017central}. Let $Y_1,\dots,Y_n$ be independent random vectors in $\R^d$, where each component is centered $\E[Y_{i,j}]=0$ and $\E[Y^2_{i,j}]<\infty$. Additionally let $X_1,\dots,X_n$ be independent random vectors, where
\begin{align*}
X_i=(X_1,\dots,X_d)^T\sim \mathcal{N}(0,\Sigma)
\end{align*}
with $\Sigma:=\E[Y_i Y_i^T]$. Assume
\begin{align*}
0<c\le \lambda_{\min}\le \lambda_{\max}\le C<\infty,
\end{align*}
where $\lambda_{\min}$ and $\lambda_{\max}$ denote the minimal and maximal eigenvalues of $\Sigma$, respectively. It is crucial that both constants $c$ and $C$ do not depend on the dimension $d$.
Further, define the normalized sums
\begin{align*}
S_n^Y:=\frac{1}{\sqrt{n}}\sum\limits_{i=1}^n Y_i \text{ and } S_n^X:=\frac{1}{\sqrt{n}}\sum\limits_{i=1}^n X_i.
\end{align*}
Next, we specify the class of $s$-sparsely convex sets as defined in \citet{chernozhukov2017central}.
\begin{definition}[$s$-sparsely convex sets]
A set $A\subset \R^d$ is called $s$-sparsely convex if there exists an integer $Q>0$ and convex sets $A_Q\subset \R^d$, $q=1,\dots,Q$, such that $A=\cap^Q_{q=1}A_q$. Additionally, the indicator function of each $A_q$, $\omega\mapsto I(\omega\in A_Q)$, depends only on $s$ elements of its argument $\omega = (\omega_1,\dots,\omega_p)$.
\end{definition}
 \noindent
\citet{chernozhukov2017central} were able to prove that, under some regularity conditions, for the class of $s$-sparsely convex sets $\mathcal{A}^{sp}(s)$ it holds
\begin{align*}
\sup_{A\in\mathcal{A}^{sp}(s)}\Big|P(S_n^Y\in A)-P(S_n^X\in A)\Big|\xrightarrow{n \to \infty} 0
\end{align*}
even if $d$ is larger than $n$. They additionally provide a result for bootstrapping in high-dimensions, enabling the construction of high-dimensional confidence regions. The standard confidence region is based on hyperrectangles.
Define
\begin{align*} 
A_\infty:=\left\{x\in \R^d\big|\|x\|_{\infty}\le c_\alpha^{(\infty)} \right\}\in \mathcal{A}^{sp}(1)
\end{align*}
as a $d$-dimensional cube with edge length of $2c_\alpha^{(\infty)}$, where
\begin{align*}
c_\alpha^{(\infty)}:= q_\alpha\left(\|X\|_{\infty}\right).
\end{align*}
Relying on bootstrap to approximate the covariance structure $\Sigma$ enables the approximation of $q_\alpha\left(\|Y\|_{\infty}\right)$ by $q_\alpha\left(\|X\|_{\infty}\right)$. From now on we will only focus on the volume of specific $s$-sparsely convex sets and omit the approximation from $Y$ to $X$. The volume $V_{A_\infty}$ of $A_\infty$ is given by
\begin{align*}
V(A_\infty)=\left(2c_\alpha^{(\infty)}\right)^d.
\end{align*}
Motivated by the well known property that the volume of the $d$-ball with fixed radius approaches zero, we will use different $s$-sparsely convex sets and analyze their behavior in large dimensions.
Let $s\in\mathbb{N}$, which is fixed and does not depend on $d$. Additionally, for simplicity assume that $\frac{d}{s}=l_s\in\mathbb{N}$ and define the corresponding index sets
\begin{align*}
J_k:=\{(k-1)\cdot s + 1,\dots, k\cdot s\},\quad k=1,\dots,l_s.
\end{align*}
Next, define, for any $p\ge 1$,
\begin{align*} 
A_p:=\left\{x\in \R^d\big|\max_{1\le k\le l_s}\|x_{J_k}\|_p \le c_\alpha^{(p)} \right\}\in \mathcal{A}^{sp}(s),
\end{align*}
which is the intersection of $l_s$ orthogonal $d$-dimensional cylinders with radius $ c_\alpha^{(p)}$, where each set only depends on $s$ components (and therefore an $s$-sparsely convex set). It can be interpreted as a crude approximation of the $d$-ball. Here  
\begin{align*}
c_\alpha^{(p)}:= q_\alpha\left(\max_{1\le k\le l_s}\|X_{J_k}\|_p\right).
\end{align*}
Since the sets $J_k$ are disjoint, it immediately follows
\begin{align*}
V(A_p)&=\left(\frac{\left(2\Gamma\left(\frac{1}{p}+1\right)c_\alpha^{(2)}\right)^s}{\Gamma\left(\frac{s}{p}+1\right)}\right)^{l_s}=\left(\frac{2\Gamma\left(\frac{1}{p}+1\right)c_\alpha^{(2)}}{\Gamma\left(\frac{s}{2}+1\right)^{\frac{1}{s}}}\right)^d,
\end{align*}
which is the volume of the $s$-ball (with respect to the $\ell_p$-norm) with radius $c_\alpha^{(p)}$ to the power of $l_s$.
To compare the volumes for a growing number of dimensions $d$, we have to consider the different size of the quantiles, since they depend on $d$. The following theorem states the main result of this note. 
\begin{theorem}\label{HDCR_th1}
For all $p\ge 2$ and $s$ large enough (the specific value only depends on the bounds of the eigenvalues), it holds
\begin{align*}
\lim_{d\to\infty}\frac{V(A_p)}{V(A_\infty)}=0.
\end{align*}
Especially, the ratio is decaying exponentially in $d$.
\end{theorem}
\noindent
Therefore, the volume of the confidence set based on $A_p$ is asymptotically negligible compared to the volume of $A_\infty$.\newpage
\noindent
\begin{proof}
At first, observe that due to Theorem 3.4 from \citet{hartigan2014bounding}, we obtain for every fixed $\alpha$ and $d$ large enough
\begin{align*}
c\sqrt{\log(d)}\le c_\alpha^{(\infty)},
\end{align*}
due to
\begin{align*}
P\left( \max_{j=1,\dots,d} |X_j|\le x\right)&\le P\left(\max_{j=1,\dots,d} X_j\le x\right)+ P\left(\min_{j=1,\dots,d} X_j\ge -x\right)\\
&=2P\left(\max_{j=1,\dots,d} X_j\le x\right).
\end{align*}
Here, the constant $c$ depends on the eigenvalues of $\Sigma$ and $\alpha$. Next, remark that 
\begin{align*}
E\left[\|X_{J_k}\|_p\right]&\le E\left[\|\Sigma_k^{\frac{1}{2}}\|_p\|\Sigma_k^{-\frac{1}{2}}X_{J_k}\|_p\right]\\
&\le s^{\frac{1}{p}-\frac{1}{2}} \|\Sigma_k^{\frac{1}{2}}\|_pE\left[ \|\Sigma_k^{-\frac{1}{2}}X_{J_k}\|_2\right]\\
&\le s^{\frac{1}{p}}  \sqrt{\lambda_{\max}}
\end{align*} for any $p\ge 2$. In the last step we used 
\begin{align*}
\|\Sigma_k^{\frac{1}{2}}\|_p\le \|\Sigma_k^{\frac{1}{2}}\|_2\le \sqrt{\lambda_{\max}}, 
\end{align*}
see e.g. \citet{goldberg1987equivalence}. We can rely on basic Gaussian concentration inequalities as in Example 5.7 from \citet{boucheron2013concentration}. It holds for all $t>0$
\begin{align*}
P\left(\|X_{J_k}\|_p-\E\left[\|X_{J_k}\|_p\right]\ge t\right)\le \exp\left(-\frac{t^2}{2\|\Sigma_k^{\frac{1}{2}}\|^2_{2,p}}\right)
\end{align*}
with
\begin{align*}
\max_{1\le k\le l_s} \|\Sigma_k^{\frac{1}{2}}\|_{2,p}&=\max_{1\le k\le l_s}\sup_{y\in\R^s:\|y\|_2=1}\|\Sigma_k^{\frac{1}{2}} y\|_p\\
&\le s^{\frac{1}{p}-\frac{1}{2}}\max_{1\le k\le l_s}\sup_{y\in\R^s:\|y\|_2=1}\|\Sigma_k^{\frac{1}{2}} y\|_2\\
&\le s^{\frac{1}{p}-\frac{1}{2}}\sqrt{\lambda_{\max}}.
\end{align*}
Therefore, for 
$$\bar{x}_p:=s^{\frac{1}{p}-\frac{1}{2}}\sqrt{2 \lambda_{\max}\log\left(\frac{d}{\alpha s}\right)}+ s^{\frac{1}{p}}  \sqrt{\lambda_{\max}}$$
we obtain
\begin{align*}
P\left(\|X_{J_k}\|_p\ge \bar{x}_p\right)&= P\left(\|X_{J_k}\|_p -E\left[\|X_{J_k}\|_p\right]\ge\bar{x}_p -E\left[\|X_{J_k}\|_p\right]\right)\\
&\le P\left(\|X_{J_k}\|_p -E\left[\|X_{J_k}\|_p\right]\ge s^{\frac{1}{p}-\frac{1}{2}}\sqrt{2\lambda_{\max}\log\left(\frac{d}{\alpha s}\right)}\right)\\
&\le  \exp\left(\log\left(\frac{\alpha s}{d}\right)\left(\frac{s^{\frac{1}{p}-\frac{1}{2}} \sqrt{\lambda_{\max}}}{\|\Sigma_k^{\frac{1}{2}}\|_{2,p}}\right)^2\right)\\
&\le  \frac{\alpha s}{d}.
\end{align*}
It follows
\begin{align*}
1 - P\left(\max_{1\le k\le l_s}\|X_{J_k}\|_p\le \bar{x}_p\right)&\le \sum_{k=1}^{l_s}\left(1-P\left(\|X_{J_k}\|_p\le \bar{x}_p\right)\right)\le \alpha.
\end{align*}
Therefore, it holds that
\begin{align*}
c_\alpha^{(p)}&\le s^{\frac{1}{p}-\frac{1}{2}}\sqrt{2 \lambda_{\max}\log\left(\frac{d}{\alpha s}\right)}+ s^{\frac{1}{p}}  \sqrt{\lambda_{\max}}
\end{align*}
As a result, we directly obtain for every fixed $\alpha$, $p\ge 2$  and $s$
\begin{align*}
c_\alpha^{(p)}&\le C s^{\frac{1}{p}-\frac{1}{2}} \sqrt{\log(d)},
\end{align*}
where the constant $C$ does not depend on $s$ as long as $d$ is large enough ($\log(d)> s$). 
If we compare the volumes of $A_\infty$ and $A_p$, it holds
\begin{align*}
\left(\frac{V(A_p)}{V(A_\infty)}\right)^{\frac{1}{d}}&=\frac{\Gamma\left(\frac{1}{p}+1\right)c_\alpha^{(p)}}{\Gamma\left(\frac{s}{p}+1\right)^{\frac{1}{s}}c_\alpha^{(\infty)}}\\
&\le \frac{\Gamma\left(\frac{1}{p}+1\right)Cs^{\frac{1}{p}-\frac{1}{2}}}{\Gamma\left(\frac{s}{p}+1\right)^{\frac{1}{s}}c}\\
&<1
\end{align*}
for $d$ large enough as long as 
\begin{align*}
\Gamma\left(\frac{s}{p}+1\right)^{\frac{1}{s}}\ge s^{\frac{1}{p}-\frac{1}{2}}\Gamma\left(\frac{1}{p}+1\right)\frac{C}{c},
\end{align*}
which will be satisfied for $s$ large enough due to the faster than exponential growth rate of the gamma function.
\end{proof}

\newpage
\section{Simulation}
This section provides a simulation study to underline our theoretical findings. Let
$$X=(X_1,\dots,X_d)\sim\mathcal{N}(0,\Sigma).$$
We consider three different correlation structures
$$\Sigma_l=\left(c_l^{|i-j|}\right)_{i,j\in \{1,\dots,d\}},\quad l=1,2,3$$
with $c_1=0$, $c_2=0.5$ and $c_3=0.9$. Observe that the corresponding eigenvalues are
bounded from above by
\begin{align*}
\|\Sigma_l\|_2\le \sqrt{\|\Sigma_l\|_1\|\Sigma_l\|_\infty}=\|\Sigma_l\|_1\le \frac{1+c_l}{1-c_l}.
\end{align*}
Following the argument from \citet{rosenblum1997hardy} (p. 62) the bound is sharp in the sense that 
\begin{align*}
\lambda_{\max}(\Sigma_l) \xrightarrow{d\to\infty} \frac{1+c_l}{1-c_l}.
\end{align*}
Since the theoretical guarantees only hold for $s$ large enough (depending on the eigenvalues), we would expect to need a larger sparsity index $s$ for a larger $c_j$.
We generate $n=10^5$ independent samples of $X$ to estimate the quantiles
$c_\alpha^{(p)}$. The number is chosen large to obtain precise estimates. Afterwards, we calculate the corresponding volume of each region and plot the ratio. Since the ratio of volumes is decaying exponentially with the dimension $d$, we plot the logarithm of the ratio for given volumes. The linear behavior in all simulations supports our theoretical results.

\begin{figure}[!ht]
\center
\includegraphics[height=0.55\textheight,width=\textwidth,keepaspectratio]{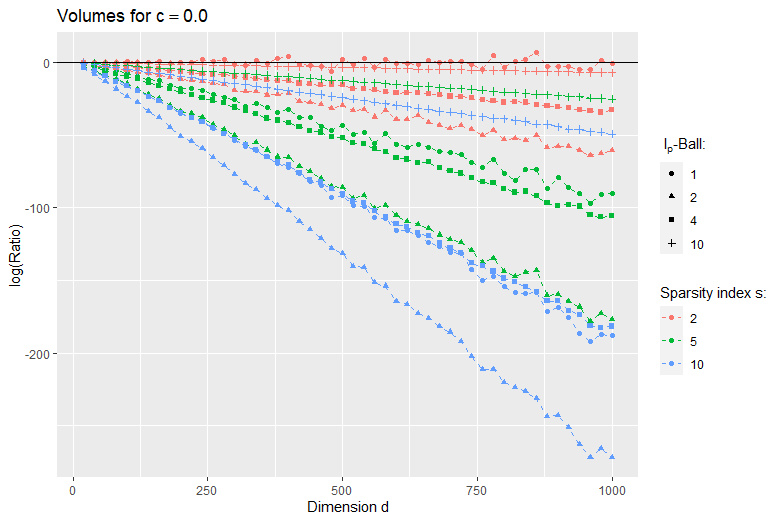}
\caption{\textmd{Simulation results for $c=0.0$.}}
\label{hdcr_c0}
\includegraphics[height=0.55\textheight,width=\textwidth,keepaspectratio]{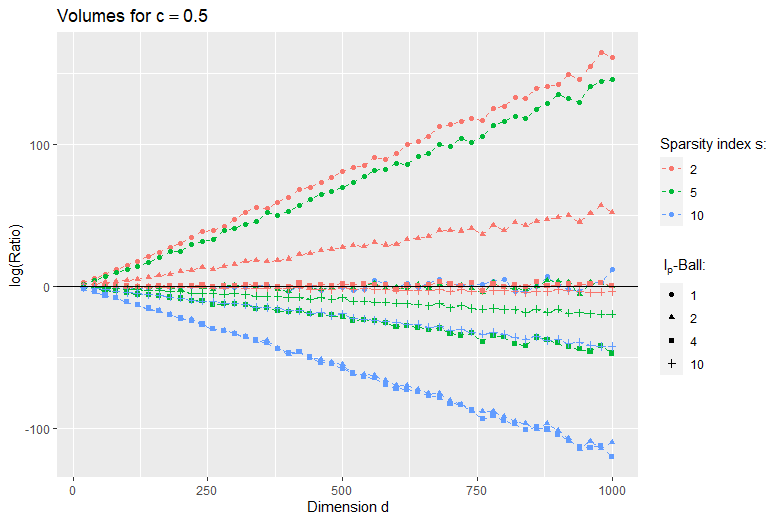}
\caption{\textmd{Simulation results for $c=0.5$.}}
\label{hdcr_c5}
\end{figure}
\
\newpage
\begin{figure}[!ht]
\center
\includegraphics[height=0.45\textheight,width=\textwidth,keepaspectratio]{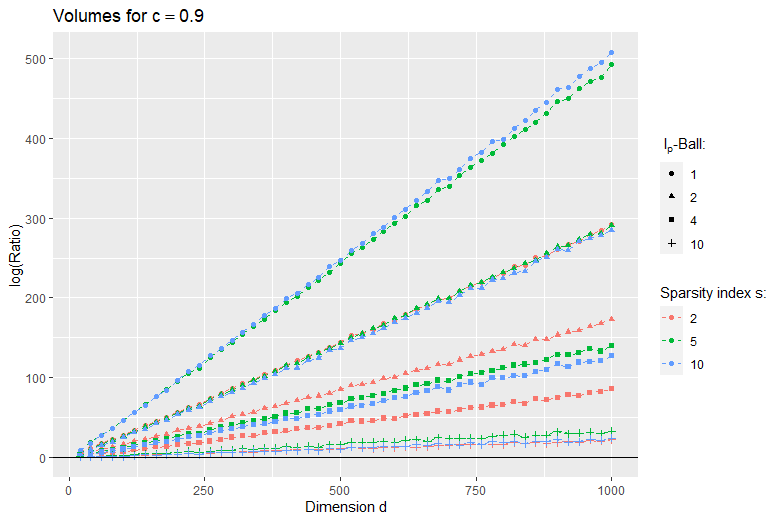}
\caption{\textmd{Simulation results for $c=0.9$.}}
\label{hdcr_c9}
\includegraphics[height=0.5\textheight,width=\textwidth,keepaspectratio]{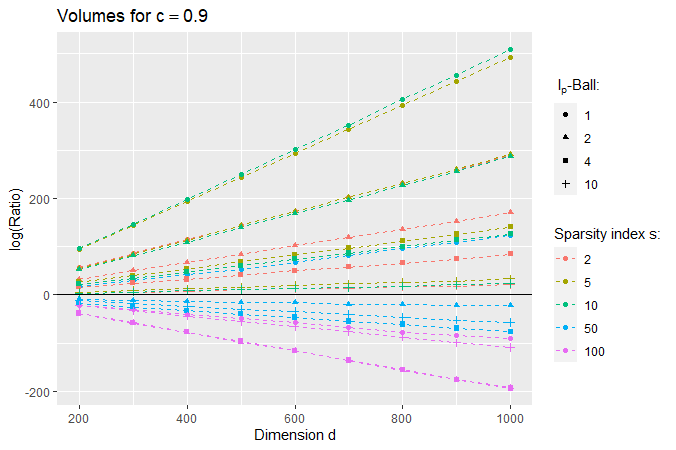}
\caption{\textmd{Additional simulation results for $c=0.9$ and larger sparsity index $s$.}}
\label{hdcr_c9_high}
\end{figure}
\noindent
In the highly correlated setting, the volume seems to be only increasing. Therefore, we add an additional plot with higher sparsity index (like Theorem \ref{HDCR_th1} would propose).
\noindent
The covariance structure has a huge effect on the corresponding quantiles (strongly positively correlated variables do not concentrate as fast as variables with weaker correlation). An simple solution to improve this problem is to permute the rows of $X$ randomly (corresponding to a randomly chosen structure of the sets $J_k$).
\begin{figure}[!htp]
\center
\includegraphics[height=0.4\textheight,width=\textwidth,keepaspectratio]{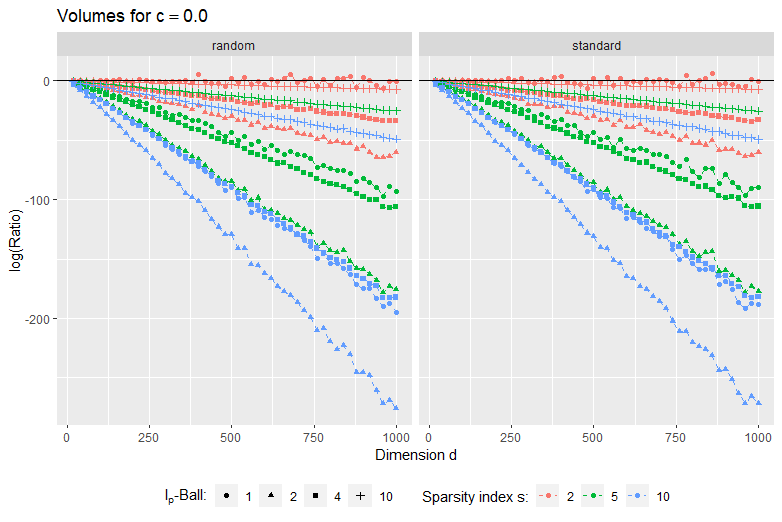}
\caption{\textmd{Simulation results for $c=0.0$ and random permutations of columns.}}
\label{hdcr_c0_comb}
\includegraphics[height=0.4\textheight,width=\textwidth,keepaspectratio]{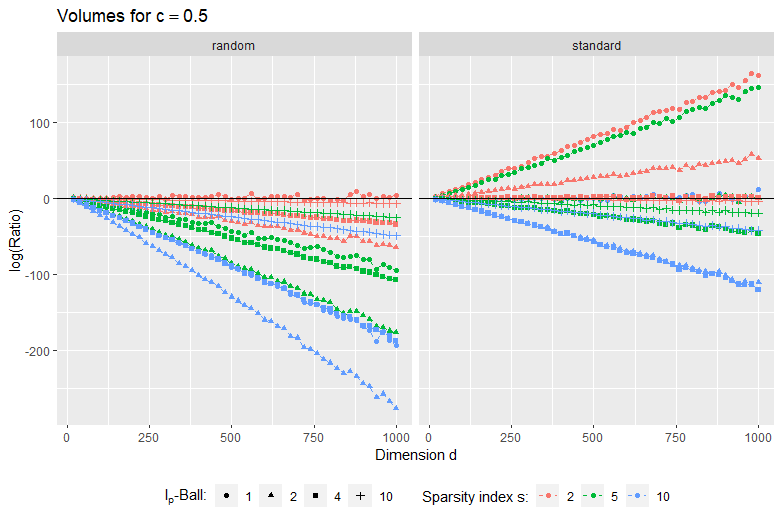}
\caption{\textmd{Simulation results for $c=0.5$ and random permutations of columns.}}
\label{hdcr_c5_comb}
\end{figure}

\begin{figure}[!htp]
\center
\includegraphics[height=0.4\textheight,width=\textwidth,keepaspectratio]{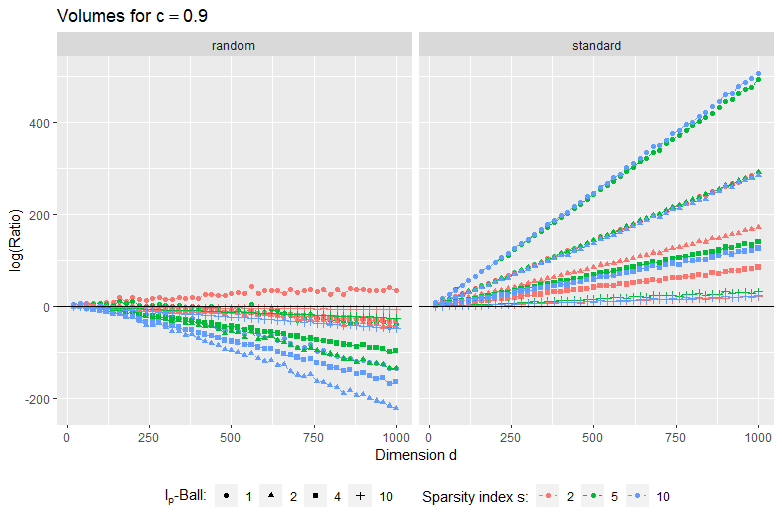}
\caption{\textmd{Simulation results for $c=0.9$ and random permutations of columns.}}
\label{hdcr_c9_comb}
\includegraphics[height=0.4\textheight,width=\textwidth,keepaspectratio]{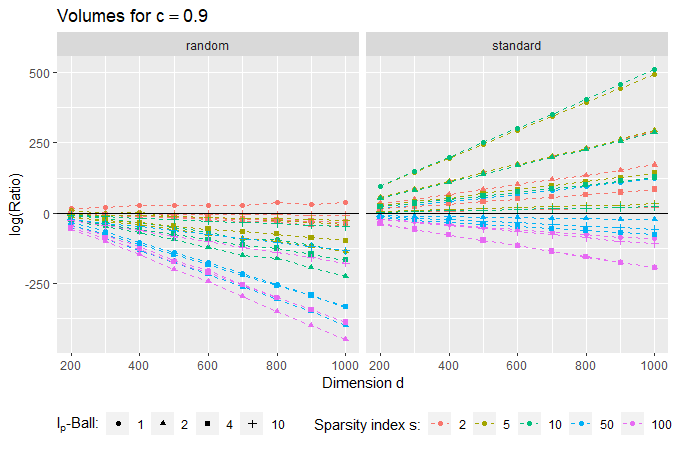}
\caption{\textmd{Additional simulation results for $c=0.9$, random permutations of columns and larger sparsity index $s$.}}
\label{hdcr_c9_high_comb}
\end{figure}
\newpage
\section{Conclusion}
In this note, we compared specific $s$-sparsely convex high-dimensional confidence regions and the corresponding hypercube with respect to their volume. Relying on Gaussian concentration inequalities, we were able to derive theoretical results demonstrating the exponential decaying ratio. In a simulation study, our theoretical results are validated as the exponential decay is clearly observable.

\newpage

\footnotesize
\pagebreak
\bibliographystyle{plainnat}
\bibliography{HDCR/literature_HDCR}

\end{document}